\documentclass[12pt]{article}
\usepackage[cp1251]{inputenc}
\usepackage[ukrainian, russian, english]{babel}

\usepackage{ amsfonts, amssymb}
\usepackage[mathscr]{eucal}

\usepackage{mathrsfs}

\begin{document}

\selectlanguage{ukrainian} \thispagestyle{empty}
 \pagestyle{myheadings}               %%%%%%%%%%%%%%%%%%%%% <--------------

\noindent {\small УДК 517.5} \vskip 1.5mm

\noindent \textbf{A.\,L.~Shidlich}\footnote{{The present investigations were supported, in part, by Grant of NAS of Ukraine for young scientists.}\vskip 1mm \hfill \textbf{\copyright\ \ A.\,L.~Shidlich, 2013}}
 {\small (Institute of Mathematics of NAS of Ukraine, Kiev)} \vskip 1.5mm

\noindent \textbf{Approximations of certain classes of functions of several variables by greedy approximants in the integral metrics}

\vskip 2mm

{\small \it \noindent We find the exact order estimates of the  approximations of the classes ${\cal F}_{q,r}^{\psi}$ of functions of several variables by greedy approximants in the integral metric. We also obtain the exact order estimates of the best $n$-term orthogonal trigonometric approximations of the classes ${\cal F}_{q,r}^{\psi}$  in the integral metric.
\hfill}

\renewcommand{\theequation}{1.\arabic{equation}}
\setcounter{equation}{0}

{\bf 1. Introduction.} Let $d$ be a fixed natural number, let $\mathbb R^d$, $\mathbb Z^d$  and $\mathbb Z^d_+$ be the sets of all ordered collections ${\bf k}:=(k_1,\ldots,k_d)$ of $d$ real,  integer and integer nonnegative numbers correspondingly. Let also $\mathbb T^d:=[0,2\pi]^d$ denote $d$-dimensional torus.

Further, let  $L_p(\mathbb T^d),~1\le p<\infty,$ be the space of all Lebesgue-measurable on $\mathbb R^d$
$2\pi$-periodic in each variable functions $f$  with finite norm
$$
||f||_{_{\scriptstyle L_p({\mathbb T}^d)}}:=\left\{\matrix{\bigg((2\pi)^{-d}\displaystyle{\int_{\mathbb T^d}}|f({\bf x})|^pd{\bf x}\bigg)^{1/p},\quad \hfill & 1\le p< \infty,\cr
\mathop{\rm ess\,sup}_{{\bf x}\in \mathbb T^d}|f({\bf x})|,\quad \hfill &  p=\infty.}\right.
$$
Set $({\bf k,x}):=k_1x_1+k_2x_2+\ldots+k_dx_d$ and for any $f\in L_1^d$, we denote the Fourier coefficients of $f$ by
$$
\widehat f({\bf k}):=(2\pi)^{-d}\int_{\mathbb T^d}f({\bf
x})e^{-i({\bf k,x})}d{\bf x},\quad{\bf k}\in\mathbb Z^d.
$$
We denote by $l_p^N$, $N=1,2,\ldots$, $0<p\le \infty$, the space ${\mathbb R}^N$ equipped with $l_p$-(quasi-)norm that is defined for ${\bf x}=\{x_k\}_{k=1}^N\in {\mathbb R}^N$ by
 $$
|{\bf x}|_p:=||{\bf x}||_{l_p}=\left\{\matrix{\bigg(\sum_{k=1}^N |x_k|^p\bigg)^{1/p},\quad 0<p<\infty,\cr \sup_{1\le i\le N}|x_i|,\quad p=\infty.}\right.
 $$
Let also  $\psi=\psi(t)$, $t\ge 1$, be a positive decreasing function, $\psi(0):=\psi(1)$ and $0<q,r\le \infty$.

In the paper, we investigate asymptotical behavior of some important approximation characteristics (in the sence of order estimates) of the classes of functions of several variables  $ {\cal F}_{q, r}^{\psi}$, defined by the following equality:
\begin{equation}\label{a2.21}%$$
{\cal F}_{q,r}^{\psi}:=\bigg\{f\in L_1({\mathbb T}^d)\ :\quad   ||\{|\widehat{f} ({\bf k})|/\psi(|{\bf k}|_{r})\}_{k\in {\mathbb Z}^d}||_{l_p({\mathbb Z}^d)}\le 1\bigg\}.
\end{equation}
If  $\psi(t)=t^{-s}$, $s\in {\mathbb N}$ and $q=1$, then ${\cal F}_{q,\infty}^{\psi}=:{\cal F}_{q,\infty}^{s}$ is a set of functions whose
$\alpha$th partial derivatives have absolutely convergent Fourier series. When $q=2$, ${\cal F}_{q,\infty}^{s}$ is equivalent  (modulo constants) to the unit ball of the Sobolev class $W^s_2$.

Approximation characteristics of the classes ${\cal F}_{q,r}^{\psi}$ for different $r\in (0,\infty]$ and for the various functions $\psi$ were investigated in the papers [1]--[8]. In particular, in [1], the authors found the exact order estimates  of the quantities of the best $n$-term trigonometric approximations of the classes ${\cal F}_{q,\infty}^{s}$, $s>0$, in the spaces $L_p(\mathbb T^d)$. In [2], V.N.~Temlyakov obtained the exact order estimates of approximations of these classes by $n$-term greedy approximants  in the spaces $L_p(\mathbb T^d)$. If $\psi(t)=R^{-t}$, $R>1$, then the exact order estimates of the quantities of the best $n$-term orthogonal trigonometric approximations in the spaces $L_p(\mathbb T^d)$, $2\le p<\infty$, of the classes ${\cal F}_{q,1}^{\psi}$ were found by V.S.~Romanyuk [3]. In the case, where $\psi(t)$ is a positive function that decreases to zero no faster than some power function, the quantities of the best $n$-term trigonometric approximations and the quantities of approximations by  $n$-term Greedy approximants of the classes ${\cal F}_{q,\infty}^{\psi}$ were studied in [4] and [5].

It should be noted that in [7; 8, Ch. XI] A.I.~Stepanets  got the exact values the best $n$-term trigonometric approximations of the classes ${\cal F}_{q,r}^{\psi}$ in the known spaces $S^p$. These results are significantly used by us in the proof and presented in Section 4.

\vskip 2mm

\renewcommand{\theequation}{2.\arabic{equation}}
\setcounter{equation}{0}

{\bf 2. Approximation characteristics.} In this section, we give the definition of  the approximation quantities for the functions of the classes ${\cal F}_{q,r}^{\psi}$, which are considered in this paper. First, for further convenience, we formulate the definition of the spaces $S^p(\mathbb T^d)$.

The space $S^p({\mathbb T}^d),$ $0<p<\infty,$ (see, for example, [8, Ch. XI]) is the space of all functions $f\in L_1(\mathbb T^d)$ such that
\begin{equation}\label{a2.1}%$$
 ||f||_{_{\scriptstyle S^{p}(\mathbb T^d)}}:=||\{|\widehat{f} ({\bf k})|\}_{k\in {\mathbb Z}^d}||_{l_p({\mathbb Z}^d)}=\bigg(\sum_{{\bf k}\in\mathbb
 Z^d}|\widehat f({\bf k})|^p\bigg)^{1/p}<\infty.
\end{equation}
The functions $f\in L_1(\mathbb T^d)$ and $g\in L_1(\mathbb T^d)$ are equivalent in the space $S^{p}(\mathbb T^d)$, when $\|f-g\|_{_{\scriptstyle S^{p}(\mathbb T^d)}}=0.$

Furter, for  $f\in L_1(\mathbb T^d)$, let $\{{\bf k}(l)\}_{k=1}^\infty=\{{\bf k}(l,f)\}_{k=1}^\infty$ denote the rearrangement of numbers ${\bf k}\in {\mathbb Z}^d$ such that
\begin{equation}\label{a2.3a1}%$$
|\widehat{f}({\bf k}(1))|\ge |\widehat{f}({\bf k}(2))|\ge \ldots.
\end{equation}
In general case, this rearrangement is not unique. In such a case, we take any rearrangement satisfying (\ref{a2.3a1}).

In the paper, the main approximation quantities for the functions $f\in {\cal F}_{q,r}^{\psi}$ are the following quantities:
\begin{equation}\label{a2.3a2}%$$
||f-G_n(f)||_{_{\scriptstyle X}}:=||f(\cdot)-\sum\limits_{l=1}^n \widehat{f}({\bf k}(l))e^{i({\bf k}(l),\cdot)}||_{_{\scriptstyle X}},
\end{equation}
\begin{equation}\label{a2.3a21}%$$
e_n^\perp(f)_{_{\scriptstyle X}}:=\inf\limits_{\gamma_n}||f(\cdot)-\sum\limits_{{\bf k}\in \gamma_n}\widehat{f}({\bf k})e^{i({\bf k},\cdot)}||_{_{\scriptstyle X}},
\end{equation}
and
\begin{equation}\label{a2.3a22}%$$
e_n(f)_{_{\scriptstyle X}}:=\inf\limits_{\gamma_n, c_{\bf k}}||f(\cdot)-\sum\limits_{{\bf k}\in \gamma_n}c_{\bf k}e^{i({\bf k},\cdot)}||_{_{\scriptstyle X}},
\end{equation}
where $X$ is one of the spaces $L_p(\mathbb T^d)$, $1\le p\le \infty$, or $S^{p}(\mathbb T^d)$, $0<p<\infty$, $\gamma_n$ is a collection of $n$ different vectors from the set ${\mathbb Z}^d$, $c_k$ are any complex numbers. Here, it is assumed that the  embedding ${\cal F}_{q,r}^{\psi}\subset X$ is true.

The quantities $(\ref{a2.3a22})$ and  $(\ref{a2.3a21})$ are respectively called  by the best $n$-term trigonometric and the best $n$-term orthogonal trigonometric approximations of the function $f$ in the space $X$. The quantity $(\ref{a2.3a2})$ is called by the approximation  of the function $f$ by greedy \ approximants in the space $X$.

For any  ${\mathfrak N}\subset X$, we set
$$%\begin{equation}\label{a2.3a24}%$$%
e_n^\perp({\mathfrak N})_{_{\scriptstyle X}}:=\sup\limits_{f\in {\mathfrak N}}e_n^\perp(f)_{_{\scriptstyle X}},
$$%\end{equation}
and
$$%\begin{equation}\label{a2.3a25}%$$
e_n({\mathfrak N})_{_{\scriptstyle X}}:=\sup\limits_{f\in {\mathfrak N}}e_n(f)_{_{\scriptstyle X}}.
$$%\end{equation}
In general case, the quantities (\ref{a2.3a2}) depend of the choice of the rearrangement satisfying (\ref{a2.3a1}). So, for the unique definition, we set
\begin{equation}\label{a2.3a25g}%$$
G_n({\mathfrak N})_{_{\scriptstyle X}}:=\sup\limits_{f\in {\mathfrak N}}\inf\limits_{\{{\bf k}(l,f)\}_{l=1}^\infty}||f(\cdot)-\sum\limits_{l=1}^n \widehat{f}({\bf k}(l,f))e^{i({\bf k}(l,f),\cdot)}||_{_{\scriptstyle X}}.
\end{equation}
In (\ref{a2.3a25g}), for any function $f\in {\mathfrak N}$, we consider the infimum on all rearrangements, satisfying (\ref{a2.3a1}), but it should be noted that results, formulated in this paper, are also true for any other rearrangements, satisfying (\ref{a2.3a1}).

Research of the quantities of the form (\ref{a2.3a2})--(\ref{a2.3a22}) goes back to the paper of S.B.~Stechkin [9]. Order estimates of these quantities on different classes of functions of one and several variables were obtained by many authors. In particular, in [10] and [11], one can be found the bibliography of papers in which the similar results were obtained.

Note that for any  $f\in L_p(\mathbb T^d)$,
\begin{equation}\label{a2.3a27}% $$
e_n(f)_{_{\scriptstyle L_p(\mathbb T^d)}}\le e_n^\perp (f)_{_{\scriptstyle L_p(\mathbb T^d)}}\le ||f-G_n(f)||_{_{\scriptstyle L_p(\mathbb T^d)}}.
\end{equation}
and by virtue of (\ref{a2.1}), for any  $f\in S^{p}(\mathbb T^d)$,
\begin{equation}\label{a2.3a26}% $$
e_n(f)_{_{\scriptstyle S^{p}(\mathbb T^d)}}=e_n^\perp (f)_{_{\scriptstyle S^{p}(\mathbb T^d)}}=||f-G_n(f)||_{_{\scriptstyle S^{p}(\mathbb T^d)}}.
\end{equation}

%%%%%%%%%%%%%%%%%%%%%%%%%%%%%%%%%%%%%%%%%%%%%%%%%%%%%%%%%%%%%%%%%%%  Основні результати %%%%%%%%%%%%%%%%%%%%%%%%%%%%%%%%%%%%%%%%%%%%%%%%%%%%%%%%%%%%%%%%%

\vskip 2mm

\renewcommand{\theequation}{3.\arabic{equation}}
\setcounter{equation}{0}

{\bf 3. Main result.} The main purpose of this work is to find the dependence of the  choice of the parameters $r$, $\psi$ and $q$ on the rate of convergence to zero, as $n\to\infty$, of the approximative characteristics of the classes ${\cal F}_{q,r}^{\psi}$.

{\bf 3.1.} As mentioned above, in the case, where $\psi(t)$ is a power function, i.e., $\psi(t)=t^{-s}$, $s>0$, for all $1\le p\le \infty$, the exact order estimates of the quantities $e_n({\cal F}_{q,\infty}^{\psi})_{_{\scriptstyle L_p(\mathbb T^d)}}$ and $G_n({\cal F}_{q,\infty}^{\psi})_{_{\scriptstyle L_p(\mathbb T^d)}}$ were obtained in [1] and [2], correspondingly. In particular, from Theorem 6.1 of [1], it follows that for all $s>0$, when $0<q\le 1$, and  for all $s> d(1-\frac 1q)$, when  $1<q<\infty$, the following relation is true\footnote{Here and in what follows, for positive sequences $\alpha(n)$ and $\beta(n)$, the expression '$a(n)\asymp b(n)$'  means that there are constants
$0<K_1<K_2$ such that for any $n\in\mathbb{N}$, $\ \alpha(n)\le K_2\beta(n)$ (in this case, we write '$\alpha(n)\ll\beta(n)$') and $\alpha(n)\ge K_1\beta(n)$ (in this case, we write  '$\alpha(n)\gg \beta(n)$').}:
 \begin{equation}\label{ab1a}%
e_n({\cal F}_{q,\infty}^{s})_{_{\scriptstyle L_p(\mathbb T^d)}}\asymp  n^{-\frac sd-\frac 1q+\frac 12}, \  \ \ \ 1\le p\le \infty,
 \end{equation}
and from Theorem 3.1 of  [2], it follows that for all $s>0$, when $0<q\le 1$, and  for all $s> d(1-\frac 1q)$, when  $1<q<\infty$,
\begin{equation}\label{ab1}% $$
G_n({\cal F}_{q,\infty}^{s})_{_{\scriptstyle L_p(\mathbb T^d)}} \asymp \left\{\matrix{ n^{-\frac sd-\frac 1q+\frac 12},\  \ \ \ 1\le p<2, \cr n^{-\frac sd-\frac 1q+1-\frac 1p},\  \ \ \ 2\,\le\, p\,<\infty.}\right.
\end{equation}

From the following Theorem 3.1, in particular, it follows that the estimates of form (\ref{ab1}) of the quantities $G_n({\cal F}_{q,\infty}^{\psi})_{_{\scriptstyle L_p(\mathbb T^d)}}$ are satisfied for a wider set of the functions $\psi$. Actually, this result was obtained in   [5], where the case $r=\infty$ was considered. For completeness and in view of the relative inaccessibility of [5] for English-speaking readers, in this paper, we give it  with a proof.

To formulate this statement, we use the following notation: let $B$ denote the set of all positive descending functions such that
\begin{equation}\label{a2.2}%$$
 \lim\limits_{t\to \infty} \psi(t)=0,
\end{equation}
and for a certain number $c>1$ and for all  $t\ge 1$, the following relation is true:
\begin{equation}\label{1.2a}%$$
1<\frac{\psi(t)}{\psi(c t)}\le K.
\end{equation}
Here and in what follows  $K$, $K_\psi,\ K_1,\  K_2,\ldots$ are positive constants which are independent of the variable $t$.

{\bf Theorem 3.1.} {\it Assume that $1\le r\le \infty$, $1\le p<\infty$, $0<q<\infty$, the function $\psi$ belongs to the set $B$ and moreover
for $0<p/(p-1)<q$ and for all\ \  $t$, larger than a certain number  $t_0$, $\psi$ is convex downwards and satisfies the condition
$$%\begin{equation}\label{1.2c1213}%$$
\frac  1{\alpha(\psi,t)}\ge K_\psi>\left\{\matrix{  d(\frac 12-\frac 1q),\ \  \ 1< p\le 2, \cr d(1-\frac 1p-\frac 1q),\  2\,{\le}\, p\,{<}\infty,}\right.
$$%\end{equation}
where
 \begin{equation}\label{1.2cm}%$$
\alpha(\psi,t):=\frac{\psi(t)}{t|\psi'(t)|}, \quad \psi'(t):=\psi'(t+).
\end{equation}
Then
\begin{equation}\label{1.2c122}%$$
G_n({\cal F}_{q,r}^{\psi})_{_{\scriptstyle L_p(\mathbb T^d)}}\asymp e_n^\perp({\cal F}_{q,r}^{\psi})_{_{\scriptstyle L_p(\mathbb T^d)}}\asymp  \left\{\matrix{  \frac {\psi(n^{1/d})}{n^{\frac 1q-\frac 12}},\ \  \ 1\le p\le 2, \cr \frac {\psi(n^{1/d})}{n^{\frac 1q+\frac 1p-1}},\ \  \  2\,{\le}\, p\,{<}\infty.}\right.
\end{equation}}

Note that the conditions on the function $\psi$ in Theorem 3.1 guarantee the embedding ${\cal F}_{q,r}^{\psi}{\subset} L_p({\mathbb T}^d)$.

Putting $r=\infty$ and $\psi(t)=t^{-s}$, $s>0$, from Theorem 3.1 we obtain the following corollary:

{\it \textbf{Corollary }} \textbf{3.1.} {\it  Assume that $1\le p<\infty$, $0<q<\infty$ and $s$ is a positive number, which for
$0<p/(p-1)<q$, satisfies the condition
 $$
 s{>}\left\{\matrix{  d(\frac 12-\frac 1q),\ \  \ 1< p\le 2, \cr d(1-\frac 1p-\frac 1q),\  2\,{\le}\, p\,{<}\infty.}\right.
 $$
Then
\begin{equation}\label{ab131}% $$
G_n({\cal F}_{q,\infty}^{\psi})_{_{\scriptstyle L_p(\mathbb T^d)}}\asymp e_n^\perp({\cal F}_{q,\infty}^{s})_{_{\scriptstyle L_p(\mathbb T^d)}}\asymp \left\{\matrix{  n^{-\frac sd-\frac 1q+\frac 12},\  \ \ \ 1\le p<2, \cr n^{-\frac sd-\frac 1q+1-\frac 1p},\  2\,{\le}\, p\,{<}\infty.
}\right.
\end{equation}}

For $1\le p<\infty$, this statement complements the result mentioned above of [2]   in the following sense:

\begin{itemize}

 \item[$\bullet$] from Corollary 3.1, in particular, it follows that in the case, where  $1<q\le \frac p{p-1}$, relation (\ref{ab1}) also holds for all $s>0$,

 \item[$\bullet$] if $0<p/(p-1)<q$ and $1<p\le 2$, then relation (\ref{ab1})  also holds for all  $s$ such that $d(\frac 12-\frac 1q)<s\le d(1-\frac 1q)$,

 \item[$\bullet$] if $0<p/(p-1)<q$ and $2<p<\infty$, then relation (\ref{ab1})  also holds for all  $s$ such that $d(1-\frac 1p-\frac 1q)<s\le d(1-\frac 1q)$.
 \end{itemize}

 Note also that if $0<q\le p/(p-1)$, then the conditions of Theorem 3.1 are satisfied, for example, for the function $\psi(t)=t^{-s}\ln^\varepsilon (t+e)$, where $s> 0$, $\varepsilon \in {\mathbb R}$, as well as for the function $\psi(t)=\ln^\varepsilon (t+e)$,  $\varepsilon<0$. If $1<p/(p-1)<q$ and $1<p\le 2$, then the conditions of Theorem 3.1 are satisfied  for the function
 $\psi(t)=t^{-s}\ln^\varepsilon (t+e)$, where $\varepsilon \in {\mathbb R}$ and $s>d(\frac 12-\frac 1q)$.  If $1<p/(p-1)<q$ and $2<p<\infty$, then the conditions of Theorem 3.1 are satisfied  for the function $\psi(t)=t^{-s}\ln^\varepsilon (t+e)$, where $\varepsilon \in {\mathbb R}$ and $s>d(1-\frac 1p-\frac 1q)$.

The proof of Theorem 3.1 will be given in Section 5.

%%%%%%%%%%%%%%%%%%%%%%%%%%%%%%%%%%%%%%%%%%%%%%%%%%%%%%%%%%%%%%%%%%%%%%%%%%%%%%%%%%%%%%%%%%%%%%%%%%%%%%%%%%%%%%%%%%%%%%%%%%%%%%%%%%%%%%%%%%%%%%%%%%%%%%%%%
\vskip 2mm

%%%%%%%%%%%%%%%%%%%%%%%%%%%%%%%%%%%%%%%%%%%%%%%%%%%%%%%%%%%%%%%%%%%%%%%%%%%%%%%%%%%%%%%%%%%%%%%%%%%%%%%%%%%%%%%%%%%%%%%%%%%%%%%%%%%%%%%%%%%%%%%%%%%%%%%%%

\renewcommand{\theequation}{4.\arabic{equation}}
\setcounter{equation}{0}

{\bf 4. Order estimates for some functionals and their applications.}

{\bf 4.1.} Let  $\Psi=\Psi(k)$, $k=1,2,\ldots,$ be a nonincreasing positive sequence such that
\begin{equation}\label{v1.2*}%$$
\lim\limits_{k\to+\infty}\Psi(k)=0.
\end{equation}
The following Lemma 4.1 is essentially used  for proving upper estimates in Theorem 3.1. This lemma gives exact order estimates for the following functionals ${H}_n(\Psi,s)$, which in the case, where  $s\in (0,1]$, are defined by the equality
\begin{equation}\label{v1.2}%$$
{H}_n(\Psi,s):= \sup\limits_{l> n}(l-n)
  \bigg(\sum_{k=1}^l \,{\Psi^{-s}(k)}\bigg)^{-\frac 1s},
 %\eqno(\ref{1.2})
 \end{equation}
and for $s\in (1,\infty)$, they are defined by the equality
 \begin{equation}\label{v1.3}%$$
{H}_n(\Psi,s):=\bigg((l^*-n)^{s'}\bigg(\sum _{j=1}^l\Psi^{-s}(j)\bigg)^{-{{s'}/s}}+\sum\limits_{j=l^*+1}^\infty
\Psi^{s'}(j)\bigg)^{1/{s'}},
\end{equation}
where   $1/s+1/{s'}=1$,
\begin{equation}\label{1.2c201}%$$
\sum\limits_{j=1}^\infty \Psi^{s'}(j)<\infty,
\end{equation}
and the number $l^*$ is given by relation
\begin{equation}\label{1.2c20}%$$
 \Psi^{-s}({l^*})\le \frac 1{l^*-n}\sum _{j=1}^{l^*}\Psi^{-s}(j)<\Psi^{-s}({l^*+1}).
\end{equation}

Note that, in the terms of similar functionals, are formulated solutions of many problems of approximation theory (see, eg, [7; 8, Ch. XI; 12--18]). Therefore, the problem of finding such estimates  is interesting in itself.

Let $\nu=\{\nu_i\}_{i=0}^\infty$ be an increasing sequence of natural numbers, $\nu_0:=1$, $V_m=\sum\limits_{k=0}^m \nu_k$, let $I_m{:=}I_m(\nu)=\bigg(V_{m-1},V_{m}\bigg]$, $m=1,2,\ldots$,  be the set of intervals. Let also $d\in {\mathbb N}$, $M_0$, $c_1$ and  $c_2$ be fixed positive numbers.

Let $S_d(M_0)=S_d(M_0,c_1,c_2)$ denote the set of all positive nonincreasing sequences $\Psi=\Psi(k)$, $k=1,2,\ldots,$ satisfying condition (\ref{v1.2*}), which  are represented as
 \begin{equation}\label{1.11}%$$
 \Psi(t)=\psi(m),\quad t\in I_m(\nu),\ \ m=1,2,\ldots,
  %\eqno(1.11)
 \end{equation}
where $\psi$ is the decreasing sequence of different values of the sequence $\Psi$ and the set of intervals $I_m(\nu)$ is such that for all $m$, greater than a certain number   $k_0$,
 \begin{equation}\label{1.12}%$$
M_0(m-c_1)^d< V_m\le M_0(m+c_2)^d.
%\eqno(1.12)
  \end{equation}

Without loss of generality, we assume that the sequences $\psi$ are restrictions of certain positive continuous 
functions $\psi(t)$ of continuous argument $t\ge 1$ on the set of natural numbers ${\mathbb N}$.

{\bf Lemma 4.1.} {\it Let  $s\in (0,\infty)$, $d\in {\mathbb N}$, the sequence  $\Psi$ belongs to the set $S_d(M_0)$ and the sequence of its different values is a restriction of a certain function $\psi(t)\in B$ on the set ${\mathbb N}$. Furthermore, in the case  $s>1$, we also assume that for all $t$, greater than a certain number $t_0$, the function $\psi$ is convex downwards and satisfies the condition
\begin{equation}\label{1.2c12}%$$
\alpha(\psi,t)\le K_\psi <{s'}/d,\ \ \frac 1s+\frac 1{s'}=1,
\end{equation}
where $\alpha(\psi,t)$ is defined in  $(\ref{1.2cm}).$ Then
 \begin{equation}\label{a1.1}%$$
{H}_n(\Psi;s)\asymp \frac
 {\psi(n^\frac 1d)}{n^{\frac 1s-1}}.
\end{equation}}

Let us note that for any $\Psi\in S_d(M_0)$, condition (\ref{1.2c12})  guarantees convergence of the series in (\ref{1.2c201}), when  $s>1$. Indeed, in this case, for all  $\tau\ge t_0$,
 $$
\frac{|\psi'(\tau)|}{\psi(\tau)}\ge \frac 1 {K_\psi \tau}> \frac {d}{{s'}\tau}.
 $$
After integrating each part of this relation in the range from $t_0$ to $t$, we obtain
$\psi(t)\ll t^{-1/K_\psi}\ll t^{-d/{s'}}$, $t>t_0$. Therefore, in view of (\ref{1.12}), we conclude that
$$
\sum\limits_{j=1}^\infty \Psi^{s'}(j)=  \sum\limits_{k=1}^\infty{\nu_k}\psi^{s'}(k)\ll \sum\limits_{k=1}^{\infty}{k^{d-1}}\psi^{s'}(k)
\ll \int\limits_1^\infty t^{d-1}\cdot t^{-{s'}/K_\psi}dt<\infty.
$$

%%%%%%%%%%%%%%%%%%%%%%%%%%%%%%%%%%%%%%%%%%%%%%%%%%%%%%%%%%%%%%%%%%%%%%%%%%%%%%%%%%%%%%%%%%%%%%%%%%%%%%

Also note that in the case, where the sequences $\Psi=\Psi(k)$ are restrictions of certain positive convex 
functions $\psi(t)$ of continuous argument $t\ge 1$ on the set ${\mathbb N}$, order estimates for the functionals $H_n(\Psi,r)$ were obtained in [17; 18] and order estimates for their integral analogues were  obtained in [16].

%%%%%%%%%%%%%%%%%%%%%%%%%%%%%%%%%%%%%%%%%%%%%%%%%%%%%%%%%%%%%%%%%%%%%%%%%%%%%%%%%%%%%%%%%%%%%%%%%%%%%%

\vskip 2mm

{\textbf{4.2. Proof of Lemma 4.1.}}  First, consider the case $s\in (0,1]$. In view of (\ref{1.11}), the functionals ${H}_n(\Psi;s)$ can be represented as
 $$
{H}_n(\Psi,s)= \sup\limits_{l> n}(l-n)
  \bigg(\sum_{j=1}^l \,\frac 1{\Psi^s(j)}\bigg)^{-\frac 1s}=
$$
\begin{equation}\label{a1.3}
= \sup\limits_{l> n}(l-n)
  \bigg(\sum_{k=1}^{k_l-1} \,\frac {\nu_k}{\psi^s(k)}+
  \frac{l-V_{k_l-1}}{\psi^s(k_l)}\bigg)^{-\frac
  1s}=:\widetilde{H}_n(\psi,s),
\end{equation}
where $k_l$ denote a number such that $l\in I_{k_l}$, i.e.,
\begin{equation}\label{f0}%$$
V_{k_{l}-1}<l\le  V_{k_{l}}.
  \end{equation}
By virtue of  (\ref{1.12}), for all $l>n\ge k_0$, we see that
    \begin{equation}\label{a1.9}
(l/M_0)^\frac 1d-c_2\le k_l< (l/M_0)^\frac 1d+c_1+1.
\end{equation}
From relation (\ref{1.12}), it follows that
 \begin{equation}\label{a1.4}
 \nu_k\asymp k^{d-1},
\end{equation}
therefore, for any $r>0$,
 \begin{equation}\label{a1.5}
 \sum_{k=1}^{l} \,\frac {\nu_k}{\psi^s(k)}\asymp \sum_{k=1}^{l} \,\frac {k^{d-1}}{\psi^s(k)}
\end{equation}

If $\psi\in B$, then for any $l=2,3,\ldots$, 
$$%  \begin{equation}\label{a1.6}
\frac {l^{d}}{\psi^s(l)} \ll \frac {(l/2)^{d}}{\psi^s(l/2)}  \ll \sum_{l/2\le k\le l} \,\frac {k^{d-1}}{\psi^s(k)}\le \sum_{k=1}^{l} \,\frac {k^{d-1}}{\psi^s(k)}\le \frac {l^{d}}{\psi^s(l)}.
$$%\end{equation}
Hence, according to (\ref{a1.5}), we get
\begin{equation}\label{a1.8}
 \sum_{k=1}^{l} \,\frac {\nu_k}{\psi^s(k)}\asymp \frac {l^{d}}{\psi^s(l)}.
\end{equation}

Further, by virtue of  (\ref{a1.9}) and the definition of the set $B$, we  see that
$$%\begin{equation}\label{a1.10}
\psi(k_l)\asymp\psi((l/M_0)^{\frac 1d})\asymp\psi(l^{\frac 1d}).
$$%\end{equation}
and in view of  (\ref{a1.8}), we conclude that
\begin{equation}\label{a1.11}
\widetilde{H}_n(\psi,s)\asymp \sup\limits_{l> n}(l-n)
  \bigg(\frac {k_l^d}{\psi^s(k_l)}\bigg)^{-\frac 1s}\asymp \sup\limits_{l> n}
  \bigg(\psi(l^{\frac 1d}) \frac {l-n}{l^{1/s}}\bigg).
\end{equation}

Since the function $\psi$ is non-increasing, then
\begin{equation}\label{a1.12}
\widetilde{H}_n(\psi,s)\ll \psi(n^{\frac 1d}) \sup\limits_{l> n}\frac {l-n}{l^{1/s}}.
\end{equation}
For $x>0$, $n\in {\mathbb N}$ and $s\in (0,1)$, the function $h(x)=h(x,s)=\displaystyle{\frac{x-n}{x^{1/s}}}$ attains its maximal
value at the point $x_0=n/(1-s)$, and
\begin{equation}\label{a1.13}
 h(x_0,s)=s\bigg(\frac{1-s}n\bigg)^\frac{1-s}s.
\end{equation}
If $s=1$, then the function $h(x)=h(x;1)$ is non-decreasing and tends to $1$ as $x$ increases. Therefore,
 \begin{equation}\label{a1.14}
 \sup\limits_{x>0} h(x;1)= \sup\limits_{x>0}
 \frac{x-n}{x}=\lim\limits_{x\to+\infty} \frac{x-n}{x}=1.
 \end{equation}
Combining  (\ref{a1.12})--(\ref{a1.14}), we obtain necessary upper estimates for the functionals ${H}_n(\Psi,s)$:
$$%\begin{equation}\label{a1.15}
{H}_n(\Psi,s)=\widetilde{H}_n(\psi,s)\ll \frac {\psi(n^\frac 1d)}{n^{\frac 1s-1}}.
$$%\end{equation}

Taking into account  (\ref{a1.11}) and the fact that $\psi\in B$, we also obtain the lower estimates
$$
{H}_n(\Psi,s){=}\widetilde{H}_n(\psi,s)\asymp  \sup\limits_{l> n}
  \bigg( \psi(l^{\frac 1d}) \frac {l-n}{l^{1/s}}\bigg)\ge
\psi((2n)^{\frac 1d}) \frac {2n-n}{(2n)^{1/s}}\asymp  \frac {\psi(n^\frac 1d)}{n^{\frac 1s-1}}.
$$%\end{equation}

Now, we consider tha case, when $s>1$. To simplify the notes, we set
\begin{equation}\label{a1.6311}
 Q_{n}(\Psi,l):={(l-n)}\bigg(\sum\limits_{j=1}^{l}\Psi^{-s}(j)\bigg)^{-1},\
 \ l\ge n,\ l\in {\mathbb N}.
 \end{equation}
Since for any $l>n$,
 $$
 Q_n(\Psi,l+1)=Q_n(\Psi,l)+\bigg(  \Psi^s(l+1)-Q_n(\Psi,l)\bigg)\Psi^{-s}(l+1)\bigg(
 \sum\limits_{i=1}^{l+1}\Psi^{-s}(i)  \bigg)^{  -1}
 $$
and
$$
 \Psi^s(l+1) = Q_n(\Psi,l+1)+\bigg(  \Psi^s(l+1)-Q_n(\Psi,l)  \bigg)\sum\limits_{j=1}^{l}\Psi^{-s}(j)\bigg(
 \sum\limits_{i=1}^{l+1}\Psi^{-s}(i)\bigg)^{  -1}     ,
 $$
 then in view of monotonicity of the function $\Psi$ and the definition of the number $l^*$ (see relation (\ref{1.2c20})), we conclude that for all $l\ge l^*$,
$$
 Q_n(\Psi,l)>Q_n(\Psi,l+1)>\Psi^s(l+1),
 $$
and for all $l\in [n,l^*)$,
$$
 Q_n(\Psi,l)\le Q_n(\Psi,l+1)\le \Psi^s(l+1).
 $$
This yields that
\begin{equation}\label{a1.63}
Q_n(\Psi,l^*)=\sup\limits_{l>n}  Q_{n}(\Psi,l).
 \end{equation}
According to (\ref{1.2c20}), we get $\Psi(l^*+1)>\Psi(l^*)$. Therefore, if the function $\Psi(t)$ is represented in the form as (\ref{1.11}), then
\begin{equation}\label{a1.64}
l^*=V_{k_{l^*}}=\sum\limits_{i=0}^{k_{l^*}}\nu_i
\end{equation}
where $k_{l^*}$ is defined in (\ref{f0}) for $l=l^*$.  Furthemore, in this case, the functionals  $H_n(\Psi,s)$, $r\in (1,\infty)$ can be represented as
 \begin{equation}\label{vq1.3}%$$
{H}_n(\Psi,s)=  \bigg((l^*-n)^{{s'}}\bigg(\sum _{k=1}^{k_{l^*}}\frac{\nu_k}{\psi^s(k)}\bigg)^{-{{s'} /r}}   +    \sum\limits_{k=k_{l^*}+1}^\infty
{\nu_k}{\psi^{s'}(k)}\bigg)^{1/{s'}}     :=\widetilde{H}_n(\psi,s),
\end{equation}
where $r\in (1,\infty)$, $1/s+1/{s'}=1$ and
\begin{equation}\label{q1.2c20}%$$
 \psi^{-s}({k_{l^*}})\le \frac 1{l^*-n}\sum _{j=1}^{k_{l^*}}\frac{\nu_k}{\psi^s(k)}<\psi^{-s}({k_{l^*}+1}).
\end{equation}
By virtue of (\ref{a1.63}), for the function
\begin{equation}\label{a1.65}
\widetilde{Q}_n(\psi,l):=(l-n) \bigg(\sum_{k=1}^{k_l-1} \,\frac {\nu_k}{\psi^s(k)}+ \frac{l-V_{k_l-1}}{\psi^s(k_l)}\bigg)^{-1},
\end{equation}
where $k_l$ is defined in (\ref{f0}),  the following relation is satisfied:
\begin{equation}\label{a1.66}
\sup\limits_{l> n}\widetilde{Q}_n(\psi,l)=\widetilde{Q}_n(\psi,l^*)=(l^*-n)\bigg(\sum _{k=1}^{k_{l^*}}\frac{\nu_k}{\psi^s(k)}\bigg)^{-1}.
\end{equation}

If the function  $\Psi$ satisfies the conditions of Lemma 4.1, then similarly to the case $s\in (0,1]$, we show that
\begin{equation}\label{a1.67a}
\sum_{k=1}^{l} \,\frac {\nu_k}{\psi^s(k)}\asymp \frac {l^d}{\psi^s(l)}
\end{equation}
and
\begin{equation}\label{a1.67}
\widetilde{Q}_n(\psi,l^*)=\sup\limits_{l> n}\widetilde{Q}_n(\psi,l)\asymp \psi^s(n^{\frac 1d}).
\end{equation}
From (\ref{q1.2c20}), for any $\psi\in B$,  we have
\begin{equation}\label{a1.68}
\psi(k_{l^*})\asymp \psi(n^{\frac 1d}).
\end{equation}

From relation (\ref{1.2c12}), it follows that for any $t>t_0$,
 $$
 \frac 1t\le K_\psi \frac{|\psi'(t)|}{\psi(t)}.
  $$
Integrating the left-hand and right-hand sides of this inequality in the range from a certain number $k_0$ to $k_{l^*}$,  $t_0<k_0<k_{l^*}$, we obtain
\begin{equation}\label{a1.74}
 \ln \frac {k_{l^*}}{k_0}\le K_\psi \ln \frac{\psi(k_0)}{\psi(k_{l^*})}.
\end{equation}
Putting $k_0=(n/M_0)^\frac 1d-c_2$, due to (\ref{a1.64}) and (\ref{a1.9}), we conclude that $k_0<k_{l^*}$. Therefore,
by virtue of definition of the set $B$,  relations (\ref{a1.74}) and (\ref{a1.68}), we get
\begin{equation}\label{a1.75}
 k_{l^*}\asymp n^\frac 1d.
\end{equation}

Since  $\psi\in B$, then from (\ref{a1.4}) it follows that for any $l\in {\mathbb N}$,
\begin{equation}\label{a1.69}
\sum\limits_{k=l+1}^\infty{\nu_k}\psi^s(k)\gg \sum\limits_{k=l+1}^{2l}{k^{d-1}}\psi^s(k)\gg l^{d}\psi^s(l).
\end{equation}
For $t>t_0$, the derivative of the function $h(t)=t^{d-1}\psi^{s'}(t)$  is of the form
 $$
h'(t)={s'}\psi^{s'}(t)t^{d-2}\bigg(\frac {d-1}{s'}-\frac {t|\psi'(t)|}{\psi(t)}\bigg).
 $$
Hence, taking into account (\ref{1.2c12}), we see that the function $h(t)$  decreases at $t>t_0$. Therefore, in view of (\ref{a1.4}), we obtain
\begin{equation}\label{a1.70}
\sum\limits_{k=l+1}^\infty{\nu_k}\psi^{s'}(k)\ll \sum\limits_{k=l+1}^{\infty}{k^{d-1}}\psi^{s'}(k)\ll \int\limits_l^\infty t^{d-1}\psi^{s'}(t)dt.
\end{equation}
By virtue of (\ref{1.11}), (\ref{a1.4}) and (\ref{1.2c201}),
\begin{equation}\label{1.2c201aa}%$$
\sum\limits_{j=1}^\infty \Psi^{s'}(j)=\sum\limits_{k=1}^\infty{\nu_k}\psi^{s'}(k)\asymp \sum\limits_{k=1}^\infty{k^{d-1}}\psi^{s'}(k)<\infty.
\end{equation}
Hence, due to monotonicity of the function $h(t)$, it follows that $\displaystyle{\frac{{k^{d-1}}\psi^{s'}(k)}{1/k}=k^d\psi^{s'}}(k)\to 0$ as $k\to \infty$.

Further, using (\ref{1.2c12}) and the method of integration by parts, we have
\begin{equation}\label{a1.71}
\int\limits_l^\infty t^{d-1}\psi^{s'}(t)dt\le K_\psi \int\limits_l^\infty t^{d}\psi^{{s'}-1}(t)|\psi'(t)|dt=
\frac {K_\psi}{s'}t^d\psi^{s'}(t) +\frac{K_\psi d}{s'}
\int\limits_l^\infty t^{d-1}\psi^{s'}(t)dt.
\end{equation}
Taking into account (\ref{1.2c12}), from the relations (\ref{a1.70})  and (\ref{a1.71}) we get the estimate
$$%\begin{equation}\label{a1.72}
\sum\limits_{k=l+1}^\infty{\nu_k}\psi^{s'}(k)\ll l^d\psi^{s'}(l),
$$%\end{equation}
that together with (\ref{a1.69})  proves the relation
\begin{equation}\label{a1.73}
\sum\limits_{k=l+1}^\infty{\nu_k}\psi^{s'}(k)\asymp l^d\psi^{s'}(l).
\end{equation}

Thus, combining the relations  (\ref{vq1.3}), (\ref{a1.67a})--(\ref{a1.68}), (\ref{a1.73})  and (\ref{a1.75}) we obtain the estimate (\ref{a1.1}), i.e.,
$$
{H}_n(\Psi,s)=\widetilde{H}_n(\psi,s)\asymp\bigg( \psi^{r{s'}}(n^\frac 1d)\cdot \bigg( {n}/{\psi^s(n^\frac 1d)}\bigg)^{(1-\frac 1s){s'}} 
  +
 n\psi^{s'}(n^\frac 1d)\bigg)^{1/{s'}}  \asymp
$$
$$%\begin{equation}\label{a1.76}%$$
\asymp\psi(n^\frac 1d)n^\frac 1{s'}\asymp \psi(n^\frac 1d)n^{1-\frac 1s}.
$$%\end{equation}
Lemma 4.1 is proved.

%%%%%%%%%%%%%%%%%%%%%%%%%%%%%%%%%%%%%%%%%%%%%%%%%%%%%%%%%%%%%%%%%%%%%%%%%%%%%%%%%%%%%%%%%%%%%%%%%%%%%%

%%%%%%%%%%%%%%%%%%%%%%%%%%%%%%%%%%%%%%%%%%%%%%%%%%%%%%%%%%%%%%%%%%%%%%%%%%%%%%%%%%%%%%%%%%%%%%%%%%%%%%%%%%%%%%%%%%%%%%%%%%%%%%%%%%%%%%

%%%%%%%%%%%%%%%%%%%%%%%%%%%%%%%%%%%%%%%%%%%%%%%%%%%%%%%%%%%%%%%%%%%%%%%%%%%%%%%%%%%%%%%%%%%%%%%%%%%%%%%%%%%%%%%%%%%%%%%%%%%%%%%%%%%%%%%%%%%%%%%%%%%%%%%%%
\vskip 2mm

{\bf 4.3.} In this section we apply Lemma 3.1 for estimates of approximative characteristics of the spaces $S^p({\mathbb T}^d)$. Approximative characteristics of the spaces $S^p({\mathbb T}^d)$ were studied by many authors (see, for example, [7; 8, Ch. XI; 21--25]). The exact values of the quantities $e_n({\cal F}_{q,r}^{\psi})_{_{\scriptstyle S^{p}({\mathbb T}^d)}}$, as well as  the exact values of the quantities   $G_n({\cal F}_{q,r}^{\psi})_{_{\scriptstyle S^{p}({\mathbb T}^d)}}$ and $e_n^\perp ({\cal F}_{q,r}^{\psi})_{_{\scriptstyle S^{p}({\mathbb T}^d)}}$ (due to (\ref{a2.3a26})),  for any  $0<p,\,q<\infty$, were obtained by A.I. Stepanets ([7; 8, Ch. XI]). In particular, from Theorem 9.1 of [8, Ch. XI], it follows that for any $ 0<q\le p<\infty$ and for any positive function  $\psi=\psi(t)$, $t\ge 1$, satisfying condition (\ref{a2.2}), for all  $ n\in {\mathbb N}$
\begin{equation}\label{a2.3a28}% $$
e_n^p({\cal F}_{q,r}^{\psi})_{_{\scriptstyle S^{p}({\mathbb T}^d)}}=\sup\limits
 _{l>n}(l-n)(\sum _{j=1}^l\bar \psi^{-q}(j))^{-\frac
 p{q}},
\end{equation}
where  $ \bar \psi =\bar \psi(j)$, $j=1,2,\ldots$, is  the decreasing rearrangement of the system of numbers $\psi(|{\bf k}|_r)$, ${\bf k}\,{\in} {\mathbb Z}^d$. If $0<p<q<\infty$ and the positive function  $\psi=\psi(t)$, $t\ge 0$, satisfies the condition
 \begin{equation}\label{a2.3a29}% $$
\sum\limits_{{\bf k}\in {\mathbb Z}^d} \psi^\frac{pq}{q-p}(|{\bf k}|_r)<\infty,
 \end{equation}
then from Theorem 9.4 of [8, Ch. XI] it follows that
\begin{equation}\label{a2.3a30}% $$
e_n^p({\cal F}_{q,r}^{\psi})_{_{\scriptstyle S^{p}({\mathbb T}^d)}}= \bigg( (l^*-n)^{\frac{q}{q-p}}\bigg(\sum\limits_{k=1}^{l^*} \bar{\psi}^{-q}(k)\bigg)^{\frac{p}{q-p}}
  +    \sum\limits_{k=l^*+1}^{\infty} \bar{\psi}^\frac{pq}{q-p}(k)\bigg)^{\frac{q-p}{q}}  ,
  \end{equation}
where  $ \bar \psi =\bar \psi(j)$, $j=1,2,\ldots$, is  the decreasing rearrangement of the system of numbers  $\psi(|{\bf k}|_r)$, ${\bf k}\,{\in} {\mathbb Z}^d$, and the number $l^*$ is defined by
$$% \begin{equation}\label{a2.3a31}% $$
   \bar{\psi}^{-q}(l^*) \leq
   \frac{1}{l^*-n}\,\sum\limits_{k=1}^{l^*} \bar{\psi}^{-q}(k) < \bar{\psi}^
   {-q}(l^*+1).
 $$% \end{equation}

Taking into account notation (\ref{v1.2}) and (\ref{v1.3}), we can write relations (\ref{a2.3a28}) and (\ref{a2.3a30}) as
$$%  \begin{equation}\label{a2.3a32}% $$%
 e_n^p({\cal F}_{q,r}^{\psi})_{_{\scriptstyle S^{p}({\mathbb T}^d)}}=H_n(\bar\psi^p,q/p),\quad 0<p,\,q<\infty.
$$% \end{equation}
Furthermore, if the number $V_m:=|\widetilde{\Delta}_{m,r}^d|$ of elements of the set
\begin{equation}\label{f300}% $$
\widetilde{\Delta}_{m,r}^d:=\{{\bf k}\in {\mathbb Z}^d:|{\bf k}|_r\le  m,\ \ \ m\in\mathbb{N}\}.
 \end{equation}
for all sufficiently large $m\in\mathbb{Z}_+$ ($m$ greater than some positive number $k_0$) satisfies the following condition:
\begin{equation}\label{a2.212}%$$
{M_0(m-c_1)^d}<V_m= |\widetilde{\Delta}_{m,r}^d|\le {M_0(m+c_2)^d},
\end{equation}
where  $M_0$, $c_1$ and  $c_2$ are certain positive constants, then the sequence  $\bar\psi=\bar\psi(j)$, $j=1,2,\ldots$, belongs to the set $S_d(M_0)=S_d(M_0,c_1,c_2)$. Thus, by virue of Lemma 4.1, we can formulate the following statement:

\textbf{  Assertion 4.1.} {\it Assume that  $0<r\le\infty$, $1\le p<\infty$, $0<q<\infty$,  condition $(\ref{a2.212})$ holds, the function      $\psi^p(\cdot)$ belongs to the set $B$ and moreover for $0<p<q$ and for all\ \  $t$, larger than a certain number  $t_0$, $\psi^p(\cdot)$ is convex downwards and satisfies the relation
\begin{equation}\label{1.2c121}%$$%
\frac 1{\alpha(\psi,t)}\ge K_\psi>d\bigg(\frac 1p-\frac 1q\bigg).
\end{equation}
Then
$$
e_n({\cal F}_{q,r}^{\psi})_{_{\scriptstyle S^{p}({\mathbb T}^d)}}\asymp
 \frac {\psi(n^\frac 1d)}{n^{\frac 1q-\frac 1p}}.
 $$}

It is clear that in the case, when $r=\infty$, condition (\ref{a2.212}) is satisfied and $M_0=vol\{{\bf k}\in {\mathbb R}^d:|{\bf k}|_{\infty}{\le}1\}{=}2^d$. If $r=1$, then $M_0=vol\{{\bf k}\in {\mathbb R}^d:|{\bf k}|_{1}\le  1\}=2^d/d!$. Therefore, in these cases, we have
 $$
e_n({\cal F}_{q,1}^{\psi})_{_{\scriptstyle S^{p}({\mathbb T}^d)}}\asymp e_n({\cal F}_{q,\infty}^{\psi})_{_{\scriptstyle S^{p}({\mathbb T}^d)}}\asymp
 \frac {\psi(n^\frac 1d)}{n^{\frac 1q-\frac 1p}}
 $$
(for $p$, $q$ and $\psi$, satisfying conditions of Assertion 4.1). Unfortunately, we do not know whether a similar relation  for other $r$ is valid. However, one can formulate the following corollary:

\textbf{  Corollary 4.1.} {\it Assume that  $1\le p<\infty$, $0<q<\infty$, the function      $\psi^p(\cdot)$ belongs to the set $B$ and moreover for $0<p<q$ and for all\ \  $t$, larger than a certain number  $t_0$, $\psi^p(\cdot)$ is convex downwards and satisfies the relation $(\ref{1.2c121})$. Then for all  $1\le r\le\infty$,
$$
e_n({\cal F}_{q,r}^{\psi})_{_{\scriptstyle S^{p}({\mathbb T}^d)}}\asymp
 \frac {\psi(n^\frac 1d)}{n^{\frac 1q-\frac 1p}}.
 $$}

Indeed, for any numbers $r\in [1,\infty]$, $0<q<\infty$ and for any positive decreasing function $\psi=\psi(t)$, $t\ge 1$,
 \begin{equation}\label{er1}%$$%
{\cal F}_{q,1}^{\psi}\subset {\cal F}_{q,r}^{\psi}\subset {\cal F}_{q,\infty}^{\psi}.
\end{equation}
Therefore, if conditions of Corollary 4.1 are satisfied, then for all  $r\in [1,\infty]$,
 $$
\frac {\psi(n^\frac 1d)}{n^{\frac 1q-\frac 1p}}\ll e_n({\cal F}_{q,1}^{\psi})_{_{\scriptstyle S^{p}({\mathbb T}^d)}}\le e_n({\cal F}_{q,r}^{\psi})_{_{\scriptstyle S^{p}({\mathbb T}^d)}}\le
e_n({\cal F}_{q,\infty}^{\psi})_{_{\scriptstyle S^{p}({\mathbb T}^d)}}\ll \frac {\psi(n^\frac 1d)}{n^{\frac 1q-\frac 1p}}.
 $$

%%%%%%%%%%%%%%%%%%%%%%%%%%%%%%%%%%%%%%%%%%%%%%%%%%%%%%%%%%%%%%%%%%%%%%%%%%%%%%%%%%%%%%%%%%%%%%%%%%%%%%

\vskip 3mm

%%%%%%%%%%%%%%%%%%%%%%%%%%%%%%%%%%%%%%%%%%%%%%%%%%%%%%%%%%%%%%%%%%%%%%%%%%%%%%%%%%%%%%%%%%%%%%%%%%%%%%%%%%%%%%%%%%%%%%%%%%%%%%%%%%%%

\renewcommand{\theequation}{5.\arabic{equation}}
\setcounter{equation}{0}

{\bf 5. Proof of Theorem 3.1.}

{\bf 5.1.} In this section we give the proof of Theorem 3.1, but first, we formulate one auxiliary lemma, which is interesting in itself.

{\bf Lemma 5.1.} {\it Assume that $2\le p<\infty$, $n\in {\mathbb N}$,  $\gamma_n=\{{\bf k}_1, {\bf k}_2,\ldots, {\bf k}_n\}$ is a collection
of $n$ vectors $k_i\in {\mathbb Z}^d$ such that $\gamma_n\subset [-cn^\frac 1d,cn^\frac 1d]^d$, where $c$ is a positive number. Then
 \begin{equation}\label{ss1}%$$
\bigg|\bigg|\sum\limits_{{\bf k}\in \gamma_n}e^{i({\bf k},{\bf \cdot})}\bigg|\bigg|_{L_p({\mathbb T}^d)}\asymp n^{1-\frac 1p}.
\end{equation}}

{\bf Proof.}  By virtue of the Hausdorff--Young theorem (see, for example,  [26, p.~16]), we get the upper estimate:
\begin{equation}\label{ss2}%$$
\bigg|\bigg|\sum\limits_{{\bf k}\in \gamma_n}e^{i({\bf k},{\bf \cdot})}\bigg|\bigg|_{L_p({\mathbb T}^d)}\ll \bigg|\bigg|\sum\limits_{{\bf k}\in \gamma_n}e^{i({\bf k},{\bf \cdot})}\bigg|\bigg|_{S^{p}({\mathbb T}^d)}=n^{1-\frac 1p}.
\end{equation}

Let us obtain the lower estimate. Based on the known trigonometric formulas, we have
$$
\bigg|\bigg|\sum\limits_{{\bf k}\in \gamma_n}e^{i({\bf k},{\bf \cdot})}\bigg|\bigg|_{L_p({\mathbb T}^d)}\asymp
\bigg(\int\limits_{{\mathbb T}^d} \bigg|\sum\limits_{{\bf k}\in \gamma_n}e^{i({\bf k},{\bf x})}\frac{e^{i({\bf 1},{\bf x})}-1}
{e^{i({\bf 1},{\bf x})}-1}\bigg|^pd{\bf x}\bigg)^\frac 1p=
\bigg(\int\limits_{{\mathbb T}^d} \bigg|\sum\limits_{{\bf k}\in \gamma_n}\frac {e^{i({\bf k+1},{\bf x})}-e^{i({\bf k},{\bf x})}}
{e^{i({\bf 1},{\bf x})}-1}\bigg|^pd{\bf x}\bigg)^\frac 1p=
$$
$$
=\Bigg(\int\limits_{{\mathbb T}^d} \bigg|\sum\limits_{{\bf k}\in \gamma_n}\frac{(\cos ((k_1+1)x_1+\ldots+(k_d+1)x_d)+i\sin ((k_1+1)x_1+\ldots+(k_d+1)x_d))}{\cos(x_1+\ldots+x_d)+i\sin (x_1+\ldots+x_d)-1}-
$$
$$
-\sum\limits_{{\bf k}\in \gamma_n}\frac{(\cos (k_1x_1+\ldots+k_dx_d)+i\sin (k_1x_1+\ldots+k_dx_d))}{\cos(x_1+\ldots+x_d)+i\sin (x_1+\ldots+x_d)-1}\bigg|^pd{\bf x}\Bigg)^\frac 1p=
$$
$$
=\Bigg(\int\limits_{{\mathbb T}^d} \bigg|\sum\limits_{{\bf k}\in \gamma_n}\frac{-2\sin \frac{x_1+\ldots+x_d}2 \sin
\frac{(2k_1+1)x_1+\ldots+(2k_d+1)x_d)}2} {\cos(x_1+\ldots+x_d)+i\sin (x_1+\ldots+x_d)-1}+
$$
$$
+2i\sum\limits_{{\bf k}\in \gamma_n}\frac{2\sin \frac{x_1+\ldots+x_d}2 \cos
\frac{(2k_1+1)x_1+\ldots+(2k_d+1)x_d)}2} {\cos(x_1+\ldots+x_d)+i\sin (x_1+\ldots+x_d)-1}\bigg|^pd{\bf x}\Bigg)^\frac 1p.
$$
Hence, using the definition of the module,  after simplifications we obtain 
$$
\bigg|\bigg|\sum\limits_{{\bf k}\in \gamma_n}e^{i({\bf k},{\bf \cdot})}\bigg|\bigg|_{L_p({\mathbb T}^d)}\asymp
\Bigg(\int\limits_{{\mathbb T}^d} \bigg(\sum\limits_{{\bf k}\in \gamma_n}\frac{\sin^2 \frac{x_1+\ldots+x_d}2} {\sin^2 \frac{x_1+\ldots+x_d}2}+
$$
$$
+\sum\limits_{{\bf k}\in \gamma_n} \sum\limits_{{\bf j}\not= {\bf k}}\frac{\sin^2 \frac{x_1+\ldots+x_d}2\cos
((k_1-j_1)x_1+\ldots+(k_d-j_d)x_d))}{\sin^2 \frac{x_1+\ldots+x_d}2}\bigg)^\frac p2 d{\bf x}\Bigg)^\frac 1p=
$$
\begin{equation}\label{ss3}%$$
=\Bigg(\int\limits_{{\mathbb T}^d} \bigg(n+\sum\limits_{{\bf k}\in \gamma_n} \sum\limits_{{\bf j}\not= {\bf k}}\cos
((k_1-j_1)x_1+\ldots+(k_d-j_d)x_d))\bigg)^\frac p2 d{\bf x}\Bigg)^\frac 1p.
\end{equation}

Further, we set
$$
l=l(\gamma_n)=\max\limits_{{\bf k}, {\bf j}\in \gamma_n}\max\limits_{l}|k_m-j_m|.
$$
For all $x\in [0,\pi/(2\beta)]$,  $\cos \beta x\ge 1-\frac {2\beta}\pi x$  and $\cos \alpha x>\cos \beta x$, where $0<\alpha<\beta$. Therefore, from relation (\ref{ss3}) we obtain
$$
\bigg|\bigg|\sum\limits_{{\bf k}\in \gamma_n}e^{i({\bf k},{\bf \cdot})}\bigg|\bigg|_{L_p({\mathbb T}^d)}\gg
\Bigg(\int\limits_{0}^{\frac{\pi}{2ld}}\ldots\int\limits_{0}^{\frac{\pi}{2ld}} \bigg(n+\sum\limits_{{\bf k}\in \gamma_n}
\sum\limits_{{\bf j}\not= {\bf k}} (1-\frac {2l}{\pi}(x_1+\ldots+x_d))\bigg)^\frac p2 dx_1\ldots dx_d\Bigg)^\frac 1p\gg
$$
$$
\gg \Bigg(\int\limits_{0}^{\frac{\pi}{2ld}}\ldots\int\limits_{0}^{\frac{\pi}{2ld}} \bigg(n^2-\frac {2l}{\pi}n(n-1)(x_1+\ldots+x_d)\bigg)^\frac p2 dx_1\ldots dx_d\Bigg)^\frac 1p\gg
$$
$$
\gg n\Bigg(\int\limits_{0}^{\frac{\pi}{2ld}}\ldots\int\limits_{0}^{\frac{\pi}{2ld}} \bigg(1-\frac {2l}{\pi}(x_1+\ldots+x_d)\bigg)^\frac p2 dx_1\ldots dx_d\Bigg)^\frac 1p=
$$
\begin{equation}\label{ss4}%$$
=n\bigg(\frac{\pi}{2l}\bigg)^\frac dp
\Bigg(\int\limits_{0}^{\frac{1}{d}}\ldots\int\limits_{0}^{\frac{1}{d}} \bigg(1-(x_1+\ldots+x_d)\bigg)^\frac p2 dx_1\ldots dx_d\Bigg)^\frac 1p\asymp n l^{-\frac dp}.
\end{equation}
Since $\gamma_n\subset [-cn^\frac 1d,cn^\frac 1d]^d$, then $l=l(\gamma_n)\ll n^\frac 1d$. Thus, indeed, the following estimate is true:
$$
\bigg|\bigg|\sum\limits_{{\bf k}\in \gamma_n}e^{i({\bf k},{\bf \cdot})}\bigg|\bigg|_{L_p({\mathbb T}^d)}\gg n^{1-\frac 1p}.
$$
Lemma is proved.

\vskip 2mm

{\bf 5.2.} Now, we can prove Theorem 3.1.  

{\it Upper estimates.} If  $2\le p<\infty$, then using  the Hausdorff--Young theorem and relation (\ref{a2.3a26}), we get
 \begin{equation}\label{f51}%$$
 \sup\limits_{f\in {\cal F}_{q,r}^{\psi}}||f-G_n(f)||_{_{\scriptstyle L_p({\mathbb T}^d)}}\ll
 \sup\limits_{f\in {\cal F}_{q,r}^{\psi}}||f-G_n(f)||_{_{\scriptstyle S^{p'}({\mathbb T}^d)}}=
 e_n({\cal F}_{q,r}^{\psi})_{_{\scriptstyle S^{p'}({\mathbb T}^d)}},
\end{equation}
where $\frac 1p+\frac 1{p'}=1$. In the case, where $1\le p\le 2$, we have
\begin{equation}\label{f52}%$$
 \sup\limits_{f\in {\cal F}_{q,r}^{\psi}}||f-G_n(f)||_{_{\scriptstyle L_p({\mathbb T}^d)}}\ll
 \sup\limits_{f\in {\cal F}_{q,r}^{\psi}}||f-G_n(f)||_{_{\scriptstyle L_2({\mathbb T}^d)}}=
 e_n({\cal F}_{q,r}^{\psi})_{_{\scriptstyle S^{2}({\mathbb T}^d)}},
\end{equation}
Thus, to obtain the required upper estimates, it is sufficient to use Corollary 4.1.

{\it Lower estimate.}  Let ${\cal T}_m$, $m\in {\mathbb N}$, denote the set of all polynomials of the form as
 $$
 T_m({\bf x})=\sum\limits_{|{\bf k}|_\infty\le m}\widehat{T}_m({\bf k})e^{i({\bf k},{\bf x})},
 $$
 and let ${\cal A}_q({\cal T}_m)$, $0<q<\infty$, denote the subset of all polynomials $T_m\in {\cal T}_m$ such that $||T||_{_{\scriptstyle S^{q}({\mathbb T}^d)}}\le 1$.
 From Theorem 5.2 of [1], it follows that for any $0<q<\infty$, $1\le p<\infty$, $m=1,2,\ldots$ and $n=((2m+1)^d-1)/2$,
 $$
 e_n({\cal A}_q({\cal T}_m))_{_{\scriptstyle L_{p}({\mathbb T}^d)}}\ge K n^{1/2-1/q}.
 $$

For a fixed $m\in {\mathbb N}$, consider the set
$$
\psi(dm){\cal A}_q({\cal T}_m)=\{T\in {\cal T}_m\ :\ ||T||_{_{\scriptstyle S^{q}({\mathbb T}^d)}}\le \psi(dm)\}.
$$
Due to monotonicity $\psi$, for any polynomial $T\in \psi(dm){\cal A}_q({\cal T}_m)$  we have
$$
\sum\limits_{{\bf k}\in {\mathbb Z}^d}\bigg|\frac{\widehat{T}(\bf k)}{\psi(|{\bf k}|_1)}\bigg|^q\le \sum\limits_{|{\bf k}|_\infty\le m}\bigg|\frac{\widehat{T}(\bf k)}{\psi(d|{\bf k}|_\infty)}\bigg|^q\le \sum\limits_{|{\bf k}|_\infty\le m}\bigg|\frac{\widehat{T}(\bf k)}{\psi(dm)}\bigg|^q\le
1
$$
Therefore, $\psi(dm){\cal A}_q({\cal T}_m)$ is contained in the set ${\cal F}_{q,1}^{\psi}$. In view of definition of the set $B$, for all $m=1,2,\ldots$ and $n=((2m+1)^d-1)/2$, we obtain
 $$
 e_n({\cal F}_{q,1}^{\psi})_{_{\scriptstyle L_{p}({\mathbb T}^d)}}\ge  e_n(\psi(dm){\cal A}_q({\cal T}_m))_{_{\scriptstyle L_{p}({\mathbb T}^d)}}\ge K \psi(dm)n^{\frac 12- \frac 1q}\ge K_1  \psi(n^\frac 1d)n^{\frac 12-\frac 1q}.
 $$
Taking into account the relations (\ref{a2.3a27}) and (\ref{er1}), monotonicity of the quantity $e_n$ and inclusion  $\psi\in B$, we see that for all $1\le p<\infty$ and all $1\le r\le \infty$,
\begin{equation}\label{rr1111}%$$
\sup\limits_{f\in {\cal F}_{q,r}^{\psi}}||f-G_n(f)||_{_{\scriptstyle L_{p}({\mathbb T}^d)}}\gg e_n^\perp({\cal F}_{q,r}^{\psi})_{_{\scriptstyle L_{p}({\mathbb T}^d)}}\gg e_n({\cal F}_{q,r}^{\psi})_{_{\scriptstyle L_{p}({\mathbb T}^d)}}\gg e_n({\cal F}_{q,1}^{\psi})_{_{\scriptstyle L_{p}({\mathbb T}^d)}}\gg \frac {\psi(n^\frac 1d)}{n^{\frac 1q-\frac 12}}.
\end{equation}

In the case, where $2< p<\infty$, for the quantities $e_n^\perp({\cal F}_{q,r}^{\psi})_{_{\scriptstyle L_{p}({\mathbb T}^d)}}$ and $\sup\limits_{f\in {\cal F}_{q,r}^{\psi}}||f-G_n(f)||_{_{\scriptstyle L_{p}({\mathbb T}^d)}}$, this estimate can be improved. For this purpose, consider the function
$$%\begin{equation}\label{rr1}%$$%
f_{1}({\bf x})=C_1(n)\sum\limits_{|{\bf k}|_1\le {[(2n/M_0)^{1/d}}]}e^{i({\bf k},{\bf x})},
$$%\end{equation}
where $M_0=2^d/d!$ and
$$%\begin{equation}\label{rr2}%$$
 C_1(n)=\bigg(\sum\limits_{|{\bf k}|_1\le {[(2n/M_0)^{1/d}}]}\psi^{-q}(|{\bf k}|_1)\bigg)^{-\frac 1q}.
$$%\end{equation}
It is obviously that $f_1\in {\cal F}_{q,1}^{\psi}$. Due to (\ref{a2.212}), the number $|{\Delta}_{m,1}^d|$ of elements of the set
 $$
{\Delta}_{m,1}^d:=\{{\bf k}\in {\mathbb Z}^d:|{\bf k}|_{1}=m,\ \ \ m\in\mathbb{N}\}.
 $$
 for all sufficiently large $m$ satisfies the condition
\begin{equation}\label{ss44}%$$
{M_0(m-c_3)^{d-1}}< |{\Delta}_{m,1}^d|=|\widetilde{\Delta}_{m,1}^d|-|\widetilde{\Delta}_{m-1,1}^d|\le {M_0(m-c_4)^{d-1}},
\end{equation}
 where $c_3$ and $c_4$ are some positive numbers. Therefore,
 $$%\begin{equation}\label{rr7}%$$
C_1^{-q}(n)=\sum\limits_{|{\bf k}|_1\le {[(2n/M_0)^{1/d}}]}\psi^{-q}(|{\bf k}|_1)\asymp \sum\limits_{k=1}^{[(2n/M_0)^{1/d}]}\frac {k^{d-1}}{\psi^{q}(k)}.
$$%\end{equation}
Since $\psi\in B$, then for any  $l=2,3,\ldots$,
$$%\begin{equation}\label{rr8}%$$
\frac{l^{d}}{\psi^{q}(l)}\ll \frac{(l/2)^{d}}{\psi^{q}(l/2)} \ll \sum\limits_{l/2\le k\le l}\frac {k^{d-1}}{\psi^{q}(k)}\le \sum\limits_{k=1}^{l}\frac {k^{d-1}}{\psi^{q}(k)}\ll \frac{l^{d}}{\psi^{q}(l)}.
$$%\end{equation}
This yields
$$%\begin{equation}\label{rr9}%$$
C_1(n)\asymp \psi([(2n/M_0)^{1/d}])/[(2n/M_0)^{1/d}]^{\frac dq}\asymp \psi(n^{\frac  1d})/n^{\frac 1q}.
$$%\end{equation}
In view of (\ref{a2.212}) and Lemma 5.1, for any collection $\gamma_n\subset {\mathbb Z}^d$ and the polynomial $\sum\limits_{{\bf k}\in \gamma_n}\widehat{f}_1({\bf k})e^{i({\bf k},\cdot)}$, we obtain
$$%\begin{equation}\label{rr10}%$$%
 \bigg|\bigg|f_1-\sum\limits_{{\bf k}\in \gamma_n} \widehat{f}_1({\bf k}) e^{i({\bf k},\cdot)}\bigg|\bigg|_{_{\scriptstyle L_{p}({\mathbb T}^d)}}
 =C_1(n) \bigg|\bigg| \mathop{\sum\limits_{|{\bf k}|_1\le [(2n/M_0)^{1/d}] }}\limits_{{\bf k}\not{\in} \gamma_n}     e^{i({\bf k},\cdot)}\bigg|\bigg|_{_{\scriptstyle L_{p}({\mathbb T}^d)}}  \gg C_1(n) n^{1-\frac 1p}\asymp \psi(n^{\frac 1d})/n^{\frac 1q+\frac 1p-1},
$$%\end{equation}
 Therefore, for all $2\le p<\infty$, the following estimates are true:
$$%\begin{equation}\label{rr11}%$$%
 \sup\limits_{f\in {\cal F}_{q,r}^{\psi}}||f-G_n(f)||_{_{\scriptstyle L_{p}({\mathbb T}^d)}}\gg e_n^\perp({\cal F}_{q,r}^{\psi})_{_{\scriptstyle L_{p}({\mathbb T}^d)}}\gg e_n^\perp({\cal F}_{q,1}^{\psi})_{_{\scriptstyle L_{p}({\mathbb T}^d)}}\gg e_n^\perp(f_1)_{_{\scriptstyle L_{p}({\mathbb T}^d)}}\gg \psi(n^{\frac 1d})/n^{\frac 1q+\frac 1p-1}.
$$%\end{equation}
Theorem 3.1. is proved.

{\it \textbf{  Remark }} \textbf{5.1.} Combining the relations (\ref{f51}) and (\ref{rr1111}), as well as the relations (\ref{f51}) and (\ref{rr1111}),
taking into account Corollary 4.1 and  relation (\ref{a2.3a27}), we conclude that in the case, where $1\le r\le \infty$, $0<q<\infty$ and the function $\psi$ satisfies conditions of Theorem 3.1,  for all $1\le p\le 2$,
$$
e_n({\cal F}_{q,r}^{\psi})_{_{\scriptstyle L_{p}({\mathbb T}^d)}}\asymp \frac {\psi(n^\frac 1d)}{n^{\frac 1q-\frac 12}},
$$
and for all $2<p<\infty$,
$$
 \frac {\psi(n^\frac 1d)}{n^{\frac 1q-\frac 12}}\ll e_n({\cal F}_{q,r}^{\psi})_{_{\scriptstyle L_{p}({\mathbb T}^d)}}\ll \frac {\psi(n^\frac 1d)}{n^{\frac 1q+\frac 1{p}-1}}.
$$

\vskip 3mm

\newpage
\footnotesize

%%%%%%%%%%%%%%%%%%%%%%%%%%%%%%%%%%%%%%%%%%%%%%%%%%%%%%%%%%%%%%%%%%%%%%%%%%%%%%%%%%%%%%%%%%%%%%%%%%%%%%%%%%%
\begin{itemize}

 \item[{[1]}]
 R.A. DeVore,  V.N. Temlyakov, Nonlinear approximation by
 trigonometric sums,  J. Fourier Anal. Appl., 2, No. 1 (1995) 29--48.
 \item[{[2]}] V.N. Temlyakov, Greedy Algorithm and $m$-Term Trigonometric Approximation,
Constr. Approx., 14, No. 4 (1998) 569--587.
 \item[{[3]}] V.S. Romanyuk, Nonlinear approximation of several variables functions with rapidly decreasing Fourier coefficients,
Proc. of International Conference "Differential Equations, Theory of functions and their applications"\ on the occasion of 70th anniversary of academician A.M. Samoilenko, Kyiv (2008), 95--96.
 \item[{[4]}] R. S. Li,  Y. P. Liu,  Asymptotic Estimations of $m$-term Approximation and Greedy Algorithm for Multiplier Function Classes Defined by Fourier Series, Chinese Journal of Engineering Mathematics, 25, No. 1  (2008), 90--96.
 \item[{[5]}] A.L. Shidlich,  Order estimations of Best $n$-term orthogonal trigonometric approximations of the functional classes
 {${\cal F}_{q,\infty}^{\psi}$} in the spaces  $L_p(\mathbb T^d)$, Zb. Pr. Inst. Mat. NAN Ukraine, 8, No. 1 (2011) 224--243.
 \item[{[6]}] R. S. Li,  Y. P. Liu,   Best $m$-term One-sided Trigonometric Approximation of Some
Function Classes Defined by a Kind of Multipliers, Acta Mathematica Sinica, English Series, 26, No. 5 (2010), 975--984.
 \item[{[7]}]
 A.I. Stepanets, Approximation characteristics of the spaces
 $S^p_\varphi$ in different metrics, Ukr. Mat. Zh., 53, No. 8 (2001) 1121--1146.
 \item[{[8]}]
 A.I. Stepanets, Methods of Approximation Theory, VSP, Leiden--Boston, 2005.
\item[{[9]}]
 S.B. Stechkin, On absolute convergence of orthogonal series, Dokl. Akad. Nauk SSSR (N.S.),  102 (1955), 37--40.

\item[{[10]}]
R. De Vore, Nonlinear approximation, Acta Numer., 7 (1998), 51--150.
 \item[{[11]}]
  A.S. Romanyuk, Best $M$-term trigonometric approximations of Besov classes of periodic functions of several variables. (Russian)  Izv. Ross. Akad. Nauk Ser. Mat.  67, No. 2  (2003), 61--100;  translation in  Izv. Math.  67,  No. 2  (2003),  265--302.
 \item[{[12]}]
 L. B. Sofman, Diameters of octahedra. (Russian)  Mat. Zametki,  5, No. 4  (1969), 429--436.
 \item[{[13]}]
L. B. Sofman,  Diameters of an infinite-dimensional octahedron. (Russian)  Vestnik Moskov. Univ. Ser. I Mat. Meh.,  28, No. 5  (1973), 54--56.
 \item[{[14]}]
 A. Pinkus, $n$-widths in  approximation theory, Springer-Verlag Berlin Heidelberg New York Tokyo, 1985.
 \item[{[15]}]
 Fang Gensun, Qian Lixin,  Approximation Characteristics for Diagonal Operators in Different Computational Settings, J.~Approx. Theory, 140, No. 2 (2006), 178--190.
 \item[{[16]}]
 A. I. Stepanets, A. L. Shidlich, Extremal problems for integrals of nonnegative functions. (Russian),  Izv. Ross. Akad. Nauk Ser. Mat.  74, No.3 (2010), 169--224;  translation in  Izv. Math.  74, No. 3  (2010),   607--660
 \item[{[17]}]
A. I. Stepanets, A. L. Shidlich, Extremal problems for integrals of non-negative
functions, Preprint no. 2007.2, Inst. Mat. NAS Ukraine, Kiev 2007. (Russian)
 \item[{[18]}]
 A.L. Shidlich, Order equalities for some functionals and their application to the estimation of the best $n$-term approximations and widths,  Ukr. Mat. Zh., 61, No. 10 (2009) 1403--1423.
  \item[{[19]}]
 A.I. Stepanets, A.L. Shidlich, On a criterion for convex functions, Dokl. Nats. Akad. Nauk Ukraine, No. 8 (2007) 31--36.
  \item[{[20]}]
 A.I. Stepanets, Some statements for convex functions,
Ukrain. Mat. Zh., 51,  No. 5 (1999) 688--702.
\item[{[21]}]
M.D. Sterlin, Exact constants in inverse theorems in the theory of approximation, Dokl. Akad. Nauk. SSSR, 202 (3) (1972) 545--547  (in Russian).
\item[{[22]}]
A. I. Stepanets, A. S. Serdyuk,  Direct and inverse theorems in the theory of the approximation of functions in the space $S^p$, Ukrain. Mat. Zh. 54, No. 1 (2002) 106--124; English transl. in Ukrainian Math. J. 54, No. 1 (2002)  126--148.
  \item[{[23]}]
A.I. Stepanets, Problems in approximation theory in linear spaces. (Russian)  Ukrain. Mat. Zh.  58, No. 1  (2006),  47--92;  translation in  Ukrainian Math. J.  58, No. 1  (2006), 54--102.
\item[{[24]}]
S. B. Vakarchuk,  Jackson-type inequalities and exact values of widths of classes of functions in the spaces $S^p$, $1\leq p<\infty$, Ukrain. Mat. Zh. 56 No. 5 (2004), 595--605; English transl. in Ukrainian Math. J., 56 No.5 (2004), 718--729.
\item[{[25]}]
D. M. Dyachenko,  On the properties of Fourier coefficients for functions of the class $H^\omega$,  Vestnik Moskov. Univ. Ser. I Mat. Mekh., No. 4 (2005), 18--25; English transl. in Moscow Univ. Math. Bull. 60 (4) (2005), 19–26.
  \item[{[26]}]
 V.N. Temlyakov, Approximation of Periodic Functions, Computational Mathematics and Analysis Series, Commack, New York, Nova
Science Publ.,  1993.

\end{itemize}

\end{document}